\documentclass[3p]{elsarticle}
\usepackage{natbib}
\usepackage[dvipsnames]{xcolor}
\usepackage{lineno}
\usepackage{latexsym,amsmath,amssymb,amsfonts,graphicx,graphics,epsfig,url,color}
\usepackage{float}
\usepackage{pdflscape}
\usepackage{afterpage}

\makeatletter
\def\ps@pprintTitle{%
 \let\@oddhead\@empty
 \let\@evenhead\@empty
 \def\@oddfoot{\centerline{\thepage}}%
 \let\@evenfoot\@oddfoot}
 
\begin{document}

\begin{frontmatter}

\title{Augmented Phase Reduction for Periodic Orbits Near a Homoclinic Bifurcation and for Relaxation Oscillators}

\author[mymainaddress]{Bharat Monga\corref{mycorrespondingauthor}}
\ead{monga@ucsb.edu}

\author[mymainaddress,mysecondaryaddress]{Jeff Moehlis}
\cortext[mycorrespondingauthor]{Corresponding author}
\ead{moehlis@engineering.ucsb.edu}

\address[mymainaddress]{Department of Mechanical Engineering, Engineering II Building, University of California Santa Barbara, Santa Barbara, CA 93106, United States}
\address[mysecondaryaddress]{Program in Dynamical Neuroscience, University of California Santa Barbara, Santa Barbara, CA 93106, United States}

\begin{abstract} 
Oscillators - dynamical systems with stable periodic orbits - arise in many systems of physical, technological, and biological interest. The standard phase reduction, a model reduction technique based on isochrons, can be unsuitable for oscillators which have a small-magnitude negative nontrivial Floquet exponent. This necessitates the use of the augmented phase reduction, a recently devised two-dimensional reduction technique based on isochrons and isostables. In this article, we calculate analytical expressions for the augmented phase reduction for two dynamically different planar systems: periodic orbits born out of a homoclinic bifurcation, and relaxation oscillators. To validate our calculations, we simulate models in these dynamic regimes, and compare their numerically computed augmented phase reduction with the derived analytical expressions. These analytical and numerical calculations help us to understand conditions for which the use of augmented phase reduction over the standard phase reduction can be advantageous.
\end{abstract}


\end{frontmatter}

\section{Introduction}

Periodic orbits are fundamentally important in dynamical systems theory, and they are intimately tied to other fundamental concepts such as bifurcations and chaos. Beyond their theoretical interest, they arise in many crucial physical, biological, and technological applications, from mechanical oscillations to electrical circuits to circadian
rhythms to neural activity. Standard phase reduction~\cite{winfree,guckenheimer1975,kuramoto,brow04}, a classical reduction technique based on isochrons~\cite{guckenheimer1975}, has been instrumental in understanding such oscillatory systems. It reduces the dimensionality of a dynamical system with a periodic orbit to a single phase variable, and captures the oscillator's phase change due to an external perturbation through the phase response curve (PRC). This can make the analysis of high dimensional systems more tractable, and their control \cite{moehlis06,dan14,zlotnik,tass,prc_exp2,tutorial} experimentally implementable; see e.g., \cite{stigen,nabi2013,snari15,zlotnik}. This is because although the whole state space dynamics of the system may not be known, the PRC can often be measured experimentally; see e.g., \cite{prc_exp1,prc_exp2}. 

Standard phase reduction is valid only in a small neighborhood of the periodic orbit. Consequently, the magnitude of the allowable perturbation is limited by the nontrivial Floquet exponents~\cite{guck83} of the periodic orbit: in systems with a small-magnitude negative nontrivial Floquet exponent, even a relatively small perturbation can lead to a trajectory which stays away from the periodic orbit, rendering the phase reduction inaccurate and phase reduction based control ineffective. This necessitates the use of a new reduction technique called augmented phase reduction~\cite{dan16}, a two-dimensional reduction based on both isochrons and isostables~\cite{isostables}. While the first dimension captures the phase of the oscillator along the periodic orbit, like the standard phase reduction, the second dimension captures the oscillator's transversal approach to the periodic orbit. This reduction ascertains the effect of an external stimulus on the oscillator's phase change through the PRC, and the change in its transversal distance to the periodic orbit through the isostable response curve (IRC). A similar reduction was derived in \cite{shirasaka} using Koopman operator techniques. We follow the reduction derived in \cite{dan16} for our analysis in this article. Control algorithms based on the augmented phase reduction are expected to be more effective~\cite{dan16,monga_aug}, as they can be designed to allow a larger stimulus without the risk of driving the oscillator too far away from the periodic orbit. We envision that IRCs can be measured experimentally just like PRCs, making the control based on the augmented phase reduction experimentally amenable as well.

In this paper, we analytically calculate the augmented phase reduction for periodic orbits of planar systems having distinct dynamics. Specifically, we derive expressions for relaxation oscillators and systems in which periodic orbits are born out of a homoclinic bifurcations. Our contribution is the analytical calculation of IRCs and the nontrivial Floquet exponent for each of these two systems. Our approach for the IRC calculation for a relaxation oscillator is in line with Izhikevich's analysis \cite{izhikevich} for the calculation of the PRC for such systems. To validate our calculations, we simulate two different models in these regimes, and compare their numerically computed augmented phase reduction with the derived analytical expressions. A similar analysis was done in \cite{tutorial} (cf.~\cite{guillamon2009} for an alternative approach), where we derived analytical expressions for $\lambda-\omega$ systems (including the normal form for a supercritical Hopf bifurcation and the normal form of a Bautin bifurcation which has a saddle-node bifurcation of limit cycles), and simple two-dimensional models undergoing SNIPER bifurcations. Those results together with the results of this paper give a useful catalog of analytical results for the augmented phase reduction for planar dynamical systems having a stable periodic orbit.

This article in organized as follows. In Section 2, we introduce standard and augmented phase reduction. In Section 3, we analytically calculate the augmented phase reduction for the two systems, and simulate two different models under the appropriate regimes to validate our calculations. Section 4 concludes the article by summarizing the derived analytical expressions, discussing their implications, and tabulating the analytical expressions.

\section{Standard and Augmented Phase Reduction}
In this section, we give background on the concepts of isochrons, isostables, and standard and augmented phase reduction. These concepts will be used to calculate the IRC expressions in Section 3.
\subsection{Standard Phase Reduction}
The standard phase reduction is a classical technique used to describe dynamics near a periodic orbit by reducing the dimensionality of a dynamical system to a single phase variable $\theta$~\cite{winfree,kuramoto}. Consider a general $n$-dimensional dynamical system given by 
\begin{equation}
\frac{d {\bf x}}{dt} = F({\bf x}), \qquad {\bf x} \in \mathbb{R}^n, \qquad (n \ge 2) .
\label{dxdt}
\end{equation}
Suppose this system has a stable periodic orbit $\gamma(t)$ with period $T$. For each point ${\bf x}^*$ in the basin of attraction of the periodic orbit, there exists a corresponding phase $\theta({\bf x}^*)$ such that \cite{guckenheimer1975,winfree,kuramoto,brow04,tutorial}
\begin{equation}
\lim_{t \rightarrow \infty} \left| {\bf x}(t) - \gamma \left( t+\frac{T}{2 \pi} \; \theta({\bf x}^*) \right) \right| = 0,
\end{equation}
where $\bf{x}(t)$ is the flow of the initial point $\bf{x}^*$ under the given vector field. The function $\theta({\bf x})$ is called the {\it asymptotic phase} of ${\bf x}$, and takes values in 
$[0, 2 \pi)$. {\it Isochrons} are level sets of this phase function. It is typical to define isochrons so that the phase of a trajectory advances linearly in time. This implies
\begin{equation}
\frac{d\theta}{dt} = \frac{2 \pi}{T} \equiv \omega
\label{theta_flow}
\end{equation}
both on and off the periodic orbit. Now consider the system
\begin{equation}
\frac{d {\bf x}}{dt} = F({\bf x})+U(t), \qquad {\bf x} \in \mathbb{R}^n,
\label{udxdt}
\end{equation}
where $U(t) \in \mathbb{R}^n$ is an external perturbation. Standard phase reduction can be used to reduce this system to a one dimensional system given by \cite{brow04}:
\begin{equation}
\dot \theta = \omega +\mathcal{Z}(\theta)^T U(t).
\label{spr}
\end{equation}
Here $\mathcal{Z}(\theta) \equiv \nabla_{\gamma(t)} \theta \in \mathbb{R}^n$ is the gradient of phase variable $\theta$ evaluated on the periodic orbit and is referred to as the {\it (infinitesimal) phase response curve (PRC)}. It quantifies the effect of an external perturbation on the phase of a periodic orbit. The PRC can be found by solving an adjoint equation numerically \cite{erme91,hopp97,brow04}:

\begin{equation}
\frac{d \nabla_{\gamma (t)} \theta}{d t} = 
- DF^T(\gamma(t))  \, \nabla_{\gamma (t)} \theta,
\label{adjoint_equation}
\end{equation}
subject to the initial condition
\begin{equation} 
\nabla_{\gamma (0)} \theta \cdot {F}(\gamma(0)) = \omega.
\label{e.IC}
\end{equation}
Since $\nabla_{\gamma (t)} \theta$ evolves in $\mathbb{R}^n$, \eqref{e.IC} supplies only one of $n$ required initial conditions; the rest arise from requiring that the solution $\nabla_{\gamma (t)} \theta$ to \eqref{adjoint_equation} be $T$-periodic. This adjoint equation can be solved numerically with the program XPP~\cite{xpp} to find the PRC $Q_{\rm XPP}$. Since XPP computes the PRC in terms of the change in time instead of the change in phase, we rescale the XPP PRC $Q_{\rm XPP}$ as
\[
\nabla_{\gamma} \theta = \omega Q_{\rm XPP}.
\]

Eq. (\ref{spr}) is valid only in a small neighborhood of the periodic orbit, and diverges from the true dynamics as one goes further away from the periodic orbit. Therefore, the amplitude of an external perturbation has to be small enough so that it does not drive the system far away from the periodic orbit where the phase reduction is not accurate. This limitation becomes even more important if the nontrivial Floquet exponent of the periodic orbit is a negative number small in magnitude \cite{monga_aug}. This limits the achievement of certain control objectives and thus necessitates the use of the augmented phase reduction.

\subsection{Augmented Phase Reduction} \label{APR}
\begin{figure*}[!t]
\begin{center}
\includegraphics[width=0.75\textwidth]{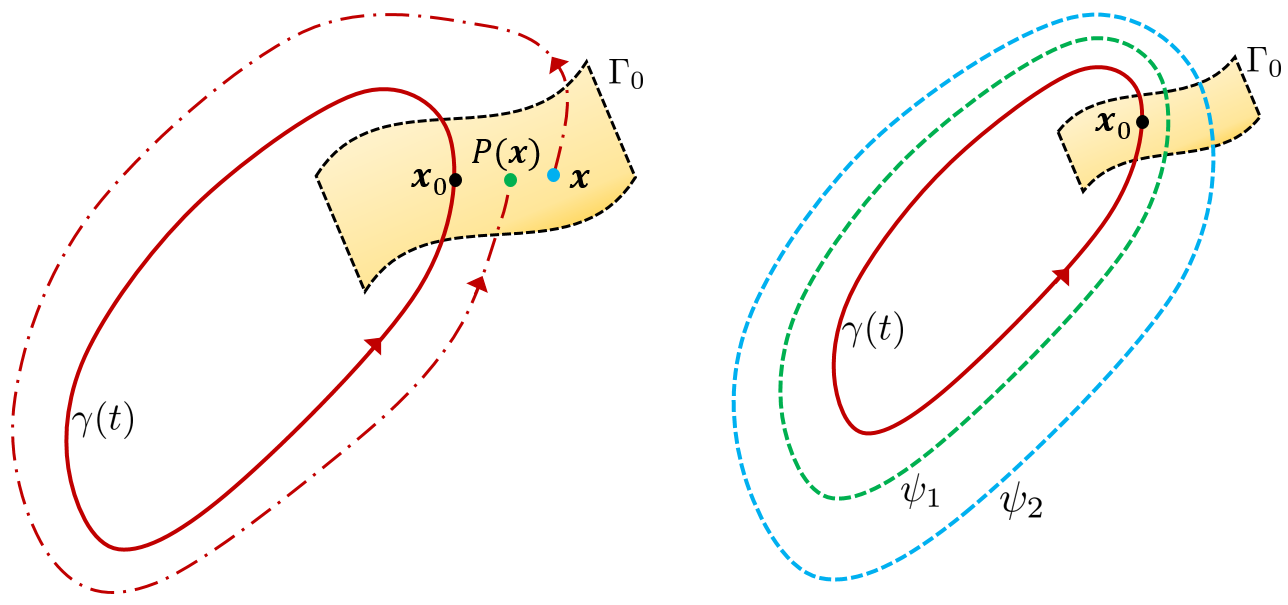}
\end{center}
\caption{Isostables for a periodic orbit. The left panel shows the Poincar\'e map $P$ on the isochron $\Gamma_0$ of the periodic orbit $\gamma (t)$. The trajectory starting from ${\bf x}$ on the isochron lands back on the isochron at $P({\bf x})$ after one period. The right panel visualizes the isostables as giving a sense of transversal distance from the periodic orbit by showing two isostable level sets $\psi_1$ and $\psi_2$.}
\label{isostable}
\end{figure*}

For systems which have a stable fixed point, it can be useful to define {\it isostables~} \cite{isostables}, which are sets of points in phase space that approach the fixed point together and are analogous to isochrons for asymptotically periodic systems. Isostables are related to the eigenfunctions of the Koopman operator \cite{isostables}. Such a notion of isostables was recently adapted for systems having a stable periodic orbit \cite{dan16}, where  isostables were defined to be the set of points that approach a periodic orbit together. They give a sense of the distance in directions transverse to the periodic orbit, visualized in the right panel of Fig. \ref{isostable}.
Standard phase reduction can be augmented with these coordinates as follows.\\
Consider a point ${\bf x_0}$ on the periodic orbit $\gamma(t)$ with the corresponding isochron $\Gamma_0$. The transient behavior of the system  (\ref{udxdt}) near ${\bf x_0}$ can be analyzed by a Poincar\'e map $P$ on $\Gamma_0$,
\begin{equation}
P:\Gamma_0 \rightarrow \Gamma_0; \qquad {\bf x} \rightarrow P({\bf x}).
\label{poin_map}
\end{equation}
This is shown in the left panel of Fig. \ref{isostable}. Here ${\bf x_0}$ is a fixed point of this map, and we can approximate $P$ in a small neighborhood of ${\bf x_0}$ as
\begin{equation}
P({\bf x})={\bf x_0}+DP({\bf x}-{\bf x_0})+{\cal} O({||{\bf x}-{\bf x_0}||}^2),
\label{lin_poin_map}
\end{equation}
where $DP=dP/dx|_{{\bf x_0}}$. Suppose $DP$ is diagonalizable with $V \in \mathbb{R}^{n\times n}$ as a matrix with columns of unit length eigenvectors $\{ v_i|i=1,\ldots,n\}$ and the associated eigenvalues $\{ \lambda_i|i=1,\ldots,n\}$ of $DP$. These eigenvalues $\lambda_i$ are the Floquet multipliers of the periodic orbit. For every nontrivial Floquet multiplier $\lambda_i$, with the corresponding eigenvector $v_i$, the set of isostable coordinates is defined as~\cite{dan16}
\begin{equation}
\psi_i({\bf x})=e_i^TV^{-1}({\bf x_\Gamma}-{\bf x_0})\exp(-\log(\lambda_i)t_\Gamma/T),
\end{equation}
where $i=1,\ldots,n-1$. Here ${\bf x_\Gamma}$ and $t_\Gamma \in [0,T)$ are defined to be the position and the time at which the trajectory first returns to the isochron $\Gamma_0$, and $e_i$ is a vector with 1 in the $i^{th}$ position and 0 elsewhere. As shown in~\cite{dan16}, we get the following equations for $\psi_i$ and its gradient $\nabla_{\gamma(t)} \psi_i$ under the flow $\dot {\bf x}=F({\bf x})$:
\begin{eqnarray}
\dot \psi_i&=&k_i \psi_i,\label{psidot}\\
\frac{d\nabla_{\gamma(t)} \psi_i}{dt}&=&\left(k_iI-DF(\gamma(t))^T\right)\nabla_{\gamma(t)}\psi_i,\label{adj1}
\end{eqnarray}
where $k_i=\log(\lambda_i)/T$ is the $i^{th}$ nontrivial Floquet exponent, $DF$ is the Jacobian of $F$, and $I$ is the identity matrix. We refer to this gradient $\nabla_{\gamma(t)}\psi_i\equiv \mathcal{I}_i(\theta)$ as the {\it isostable response curve (IRC)}. Its $T$-periodicity along with the normalization condition $\nabla_{{\bf x_0}}\psi_i \cdot v_i=1$ gives a unique IRC. It gives a measure of the effect of a control input in driving the trajectory away from the periodic orbit. The $n$-dimensional system (given by eq. (\ref{udxdt})) can be realized as~\cite{dan16}
\begin{eqnarray}
\dot \theta&=&\omega+\mathcal{Z}^T(\theta) U(t),\label{aug1}\\
\dot \psi_i&=&k_i\psi_i+\mathcal{I}_i^T(\theta) U(t), \qquad \text{for} \ i=1,\ldots,n-1.\label{aug2}
\end{eqnarray}
We refer to this reduction as the augmented phase reduction. Here, the phase variable $\theta$ indicates the position of the trajectory along the periodic orbit, and the isostable coordinate $\psi_i$ gives information about transversal distance from the periodic orbit along the $i^{th}$ eigenvector $v_i$. It is evident from (\ref{aug1}, \ref{aug2}) that an external perturbation affects the oscillator's phase through the PRC, and its transversal distance to the periodic orbit through the IRC. In practice, isostable coordinates with nontrivial Floquet multiplier close to 0 can be ignored as perturbations in those directions are nullified quickly under the evolution of the vector field. If all isostable coordinates are ignored, the augmented phase reduction reduces to the standard phase reduction. In this paper, the models that we calculate the augmented phase reduction for are two-dimensional, so there is only one isostable coordinate.  We thus write the adjoint equation as 
\begin{equation}
\frac{d\nabla_{\gamma(t)} \psi}{dt}=\left(kI-DF(\gamma(t))^T\right)\nabla_{\gamma(t)}\psi,
\label{adj}
\end{equation}
and the augmented phase reduction as
\begin{eqnarray}
\dot \theta&=&\omega+\mathcal{Z}^T(\theta)  U(t),\label{thetadot}\\
\dot \psi&=&k\psi+\mathcal{I}^T(\theta) U(t) \label{psi_dot}.
\end{eqnarray}
We have removed the subscript for $\psi$ and $k$, as we only have one isostable coordinate. The eigenvector $v$ is then the unit vector along the one-dimensional isochron $\Gamma_0$. The nontrivial Floquet exponent $k$ can then be computed from the divergence of the vector field as \cite{paul}
\begin{equation}
k=\frac{\int_0^T{ \nabla \cdot F(\gamma(t))dt}}{T}.
\label{floq_exp}
\end{equation}

\section{Analytical and Numerical Computation of the Augmented Phase Reduction} \label{section_3}
\begin{figure*}[!t]
\begin{center}
\includegraphics[width=0.75\textwidth]{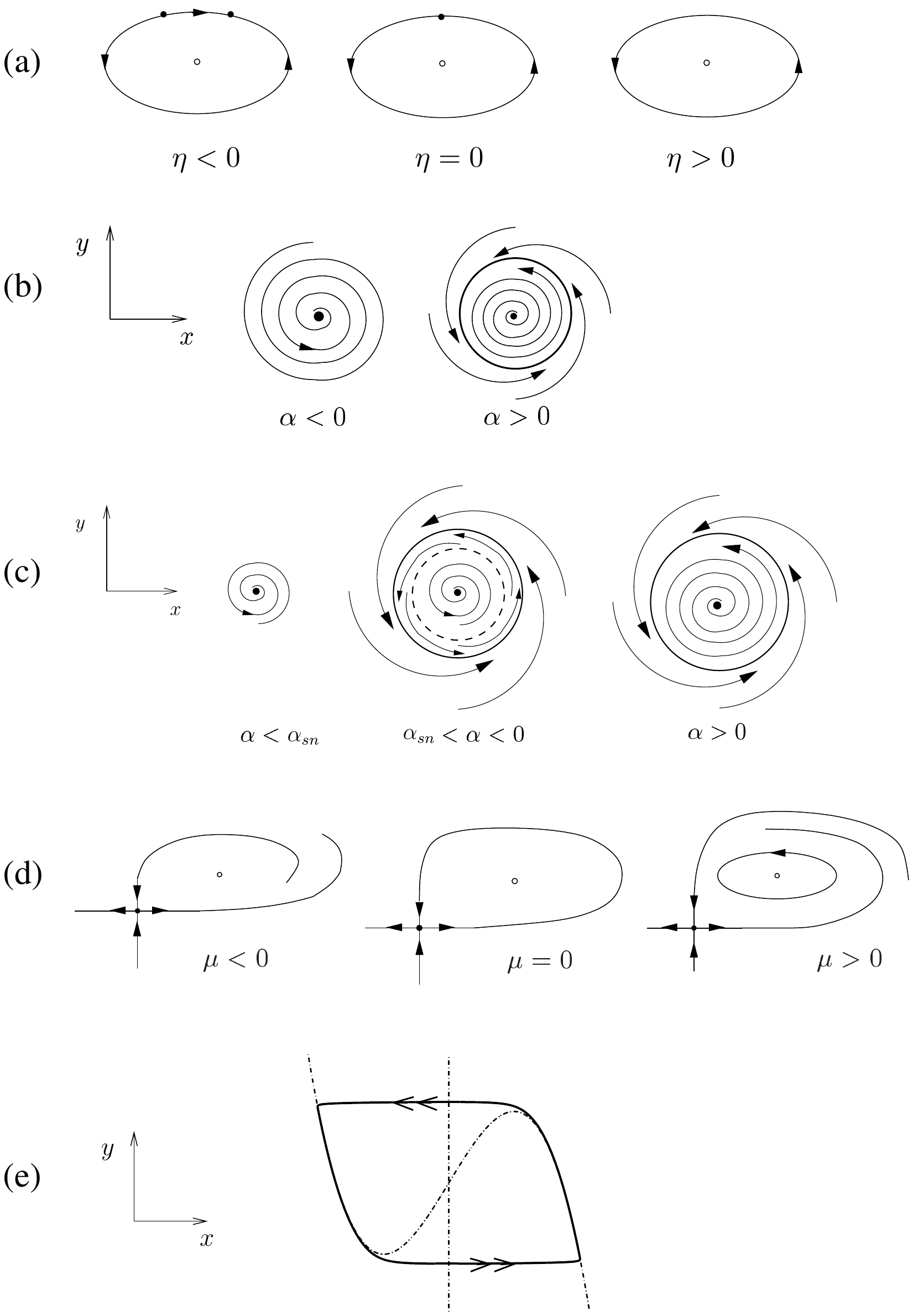}
\end{center}
\caption{(a) SNIPER bifurcation: Two fixed points die in a saddle-node bifurcation at $\eta = 1$, giving a periodic orbit for $\eta > 1$, assumed to be stable. (b) Supercritical Hopf bifurcation: A fixed point loses stability as $\alpha$ increases through zero, giving a stable periodic orbit (closed curve). (c) Bautin bifurcation: There is a subcritical Hopf bifurcation at $\alpha = 0$, and a saddle-node bifurcation of periodic orbits at $\alpha =\alpha_{sn}$. Both a stable (solid closed curve) and unstable (dashed closed curve) periodic orbit exist for $\alpha_{sn} < \alpha < 0$. The fixed point is stable (resp., unstable) for $\alpha < 0$ (resp., $\alpha > 0$). (d) Homoclinic bifurcation: A homoclinic orbit exists at $\mu = 0$, giving rise to a stable periodic orbit for $\mu > 0$. (e) A relaxation oscillator (solid closed curve) is shown with its nullclines (dashed curves).}
\label{bifurcations}
\end{figure*}

Bifurcation theory~\cite{guck83,kuznetsov} identifies four codimension one bifurcations which give birth to a stable limit cycle for generic families of vector fields: a supercritical Hopf bifurcation, a saddle-node bifurcation of limit cycles, a SNIPER bifurcation (saddle-node bifurcation of fixed points on a periodic orbit, also called a SNIC bifurcation), and a homoclinic bifurcation. These bifurcations are illustrated in Fig. \ref{bifurcations}. Fig. \ref{bifurcations} also shows a relaxation oscillator, where dynamics for one of the variables is considerably faster than the other. 

In this section, we analyze planar dynamical systems which have a stable limit cycle which arises from a homoclinic bifurcation and relaxation oscillators with fast-slow dynamics. We derive analytical expressions of the augmented phase reduction (\ref{thetadot}, \ref{psi_dot}) for these systems.

To validate our calculations, we simulate two different models whose dynamics are expected to be captured by the aforementioned planar systems. We compare their numerically computed IRCs with the derived analytical expressions. In the numerical computation of the IRCs for the planar systems, we directly calculate the nontrivial Floquet exponent $k$ as the mean of the divergence of vector field along the periodic orbit according to (\ref{floq_exp}). On the other hand, for higher dimensional models, we first compute PRC using the software XPP \cite{xpp}, then choose an arbitrary point on the periodic orbit as $\theta=0$, and approximate the isochron as a vector orthogonal to the PRC at that point. To compute the Jacobian $DF$, we compute ${\bf x_\Gamma}$ for a number of initial conditions ${\bf x_0}$ spread out on the isochron. Eigenvector decomposition of $DF$ gives us the Floquet multipliers of the periodic orbit and thus $k$. After obtaining $k$, we use Newton iteration to obtain the IRC as the periodic solution to eq. (\ref{adj}). Note that the higher dimensional systems we consider for numerical simulation in this section have only one negative small magnitude nontrivial Floquet exponent, so the reduction given by (\ref{thetadot},\ref{psi_dot}) still applies.

\subsection{Homoclinic bifurcation}
\begin{figure}[!t]
\begin{center}
\includegraphics[width=0.7\textwidth]{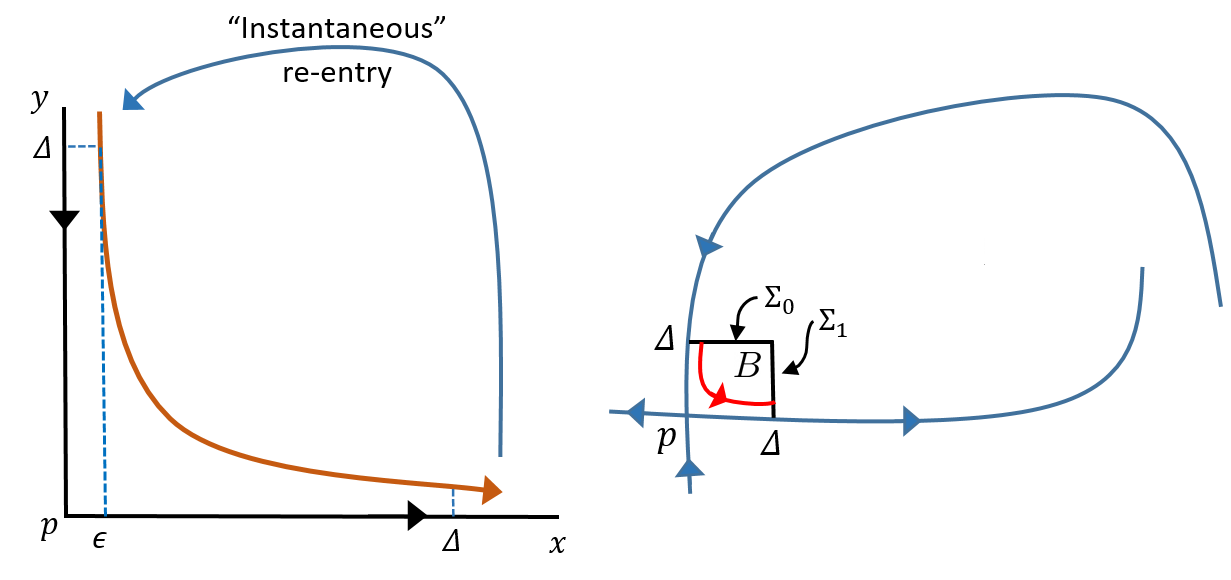}
\end{center}
\caption{Trajectory near a homoclinic bifurcation. The left panel shows the trajectory near the saddle point. The right panel shows the Poincar\'e sections used in the analysis.}
\label{box}
\end{figure}

For a homoclinic bifurcation \cite{guck83,kuznetsov}, a periodic orbit is born out of a homoclinic orbit to a hyperbolic saddle point $p$ upon varying a parameter $\mu$. If a homoclinic orbit exists for $\mu=0$, then there will be a periodic orbit for, say, $\mu > 0$, but not for $\mu<0$, as shown in Fig. \ref{bifurcations}(d). We assume that the magnitude of the unstable eigenvalue $\lambda_u$ of the saddle point is smaller than the stable eigenvalue $\lambda_s$, resulting in a stable periodic orbit \cite{guck83}. For $\mu$ close to zero, the periodic solution spends most of its time near the saddle point $p$, where the vector field can be approximated by its linearization. It can be written in diagonal form as
\begin{eqnarray}
\dot x &=& \lambda _u x,\label{lin1}\\
\dot y &=& \lambda_s y,\label{lin2}
\end{eqnarray}
where $\lambda_u >0$, and $\lambda _s < 0$. As in \cite{brow04}, we consider a box $B = [0,\Delta] \times [0,\Delta] \equiv \Sigma_0 \times \Sigma_1$ that encloses the periodic orbit for most of its time period, and within which eq. (\ref{lin1}, \ref{lin2}) are accurate. This is shown in the left panel of Fig. \ref{box}. 
We do not model the periodic orbit outside $B$, but assume that trajectory re-enters the box after a time $\delta T$ at a distance $\epsilon$ from the $y$ axis, where $\epsilon$ varies with the bifurcation parameter $\mu$. The time taken for the trajectory to traverse $B$ can be found as \cite{brow04}
\begin{equation}
\tau(\epsilon)=\frac{1}{\lambda _u}\log \left( \frac{\Delta}{\epsilon}\right).
\end{equation}
Thus the time period of the periodic orbit is given as $\tau(\epsilon)+\delta T$. As $\mu$ decreases towards zero, the periodic orbit approaches $p$, resulting in $\epsilon$ approaching 0. Near the bifurcation, $\delta T \ll \tau(\epsilon)$, so $T\approx \tau(\epsilon)$. We approximate the trajectory as spending all its time within the box $B$, and re-injecting into the box instantaneously. Thus we set $\theta = 0$ at the point where trajectory enters $B$, and $\theta = 2\pi$ where trajectory exists $B$. To find the PRC, we solve the adjoint equation in $B$ to get 
\begin{eqnarray}
\mathcal{Z}(\theta) = \mathcal{Z}_{x_0}e^{-\lambda_u t}\hat x+ \mathcal{Z}_{y_0} e^{-\lambda_s t}\hat y,
\end{eqnarray}
subject to the initial condition (eq. (\ref{e.IC}))
\begin{equation}
 \mathcal{Z}_{x_0}\lambda_u \epsilon + \mathcal{Z}_{y_0} \lambda_s \Delta = \frac{2\pi \lambda_u}{\log \left( \frac{\Delta}{\epsilon}\right)}.
\end{equation}
As $\mu \rightarrow 0$, $\epsilon \rightarrow 0$, thus the first term in the left hand side and the right hand side term in above equation go to zero. Thus we get $\mathcal{Z}_{y_0} \approx 0$ near the bifurcation point, and the PRC is only significant in the $x$-direction. Since the isochrons are orthogonal to the PRC on the limit cycle, the eigenvector $v \approx 0~\hat x + 1~\hat y$. We will use this information for the normalization condition of the IRC later. Since the trajectory spends most of its time inside the box $B$, we get $k=\lambda_s + \lambda_u$ by the mean of the divergence of the linear vector field inside $B$. We will also prove this by the Poincar\'e analysis below.

Consider the Poincar\'e maps
\begin{eqnarray}
P=&P_2 \circ P_1&:\Sigma_0 \rightarrow \Sigma_0, \qquad \text{where}\\
&P_1&:\Sigma_0 \rightarrow \Sigma_1; \qquad (x,\Delta) \rightarrow (\Delta, \Delta e^{\lambda_s T}),\\
&P_2&:\Sigma_1 \rightarrow \Sigma_0; \qquad (\Delta,y) \rightarrow (x, \Delta).
\end{eqnarray}
The Poincar\'e sections $\Sigma_0$ and $\Sigma_1$ are shown in the right panel of Fig. \ref{box}. Following the analysis in Chapter 10 of \cite{wiggins}, we get the Poincar\'e map $P$ as
\begin{equation}
P:\Sigma_0 \rightarrow \Sigma_0, \qquad (x,\Delta) \rightarrow (Ax^{-\frac{\lambda_s}{\lambda_u}}+\mu,\Delta),
\end{equation}
where $A$ is a positive constant, and $\mu$ is the bifurcation parameter. This gives the nontrivial Floquet multiplier of the periodic orbit as 
\begin{equation}
\lambda=A' \epsilon^{-\frac{\lambda_s}{\lambda_u}-1},
\end{equation}
where $A'=-A\lambda_s/\lambda_u$. From this equation, it is easy to see that $\lambda \rightarrow 0$ as $\epsilon \rightarrow 0$. Note that although the isochrons in the box $B$ may not be horizontal, we have calculated the nontrivial Floquet multiplier for a horizontal section, as that is more convenient; the value of the nontrivial Floquet multiplier is independent of the Poincar\'e section \cite{wiggins}. $k$ can be found as
\begin{equation}
k=\frac{\log\left(A' \epsilon^{-\frac{\lambda_s}{\lambda_u}-1}\right)}{T}.
\end{equation}
Near the bifurcation, this can be written as
\begin{equation}
k=\underset{\epsilon\rightarrow 0}{\lim} \frac{\log\left(A' \epsilon^{-\frac{\lambda_s}{\lambda_u}-1}\right)}{\frac{1}{\lambda _u}\log \left( \frac{\Delta}{\epsilon}\right)}.
\end{equation}
Since both the numerator and denominator approach plus or minus infinity as $\epsilon\rightarrow 0$, the limit can be solved by L'Hospital's rule as
\begin{eqnarray*}
k&=&\underset{\epsilon\rightarrow 0}{\lim} \left(\frac{\lambda _u\Delta\epsilon^{-1}}{A' \epsilon^{-\frac{\lambda_s}{\lambda_u}-1}}\right)\left(\frac{A' \left(\frac{\lambda_s}{\lambda_u}+1\right)\epsilon^{-\frac{\lambda_s}{\lambda_u}-2}}{\Delta\epsilon^{-2}}\right)\\&=&\lambda_s + \lambda_u.
\end{eqnarray*}
With this, we get the following adjoint equation for the IRC:
\begin{eqnarray}
\mathcal{\dot I}_x &=& \lambda_s\mathcal{I}_x ,\\
\mathcal{\dot I}_{y} &=& \lambda_u\mathcal{I}_y,\\
\Rightarrow  &\mathcal{I}_{x}& =\mathcal{I}_{x_0} e^{\lambda_s t},\label{ab1}\\
&\mathcal{I}_{y} &=\mathcal{I}_{y_0}e^{\lambda_u t}.\label{ab2}
\end{eqnarray}
The normalization condition $\mathcal{I}_{x_0,y_0} \ . \ v=1$ gives the IRC as
\begin{eqnarray}
 \mathcal{I}_{x,y} &=& \mathcal{I}_{x_0} e^{\frac{\lambda_s \theta}{\omega}}~\hat x + e^{\frac{\lambda_u \theta}{\omega}}~\hat y.\label{homo_irc}
\end{eqnarray}
Here $\mathcal{I}_{x_0}$ remains indeterminate as we do not model the dynamics outside $B$. The $x$ component of the IRC decreases at an exponential rate, while the $y$ component increases at an exponential rate inside the box $B$. We do not implement the condition of $T$-periodicity on eq. (\ref{ab1}, \ref{ab2}), as the calculated expressions of the IRC are valid only in the box $B$. We expect the IRC to jump back to its initial value as the trajectory re-enters the box. 

\begin{figure*}[!t]
\begin{center}
\includegraphics[width=\textwidth]{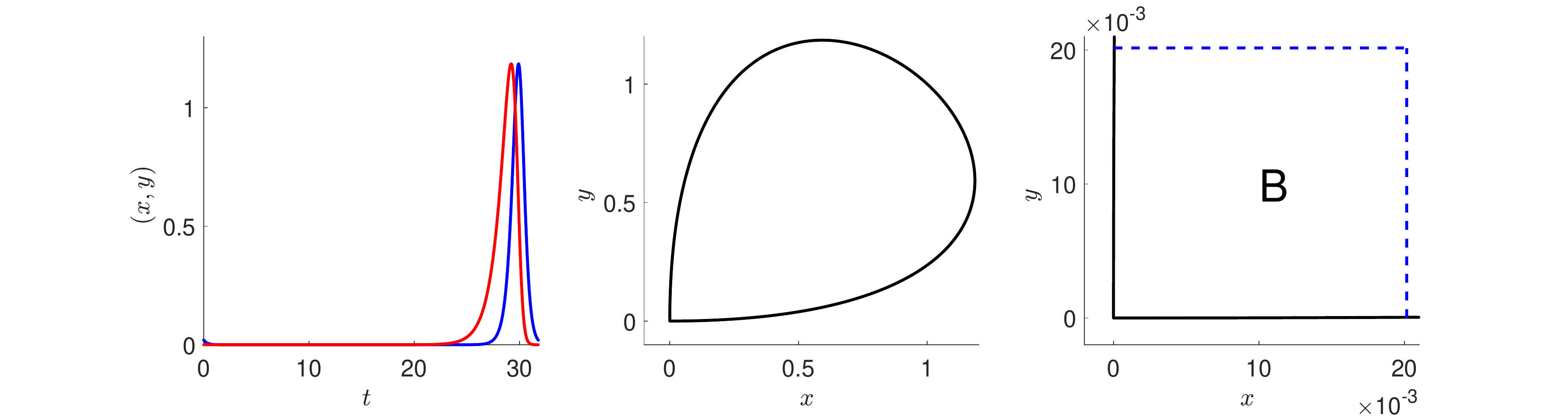}
\end{center}
\caption{Periodic orbit near homoclinic bifurcation with parameters $\mu=1 \times 10^{-13}, \ a=-1$, and $b=2$. The left (resp., middle) panel shows the time series (resp., orbit). The blue and the red lines show the $x$ and $y$ component of trajectories respectively. The right panel shows the box $B$.}
\label{simple_homoclinic}
\end{figure*}

\subsubsection{A simple model for homoclinic bifurcation}
We use  a 2-dimensional model derived from \cite{Sandstede} to validate our result:
\begin{eqnarray*}
\dot x&=&(a+b-0.5\mu)x-0.5\mu y\\
&\ & \ -(a/4+3b/8)(x+y)^2-3a/8(x^2-y^2), \ \ \\
\dot y&=&0.5\mu x+(a-b+0.5\mu)y\\
&\ & \ +(-a/4+3b/8)(x+y)^2+3a/8(x^2-y^2). \ \ \ \ \ 
\end{eqnarray*}
This system undergoes a homoclinic bifurcation at $\mu=0$, and has a stable periodic orbit for $\mu>0, \ a<0<b,$ and $|b|>|a|$. With parameters $\mu=1 \times 10^{-13}, \ a=-1$, and $b=2$, we get a stable periodic orbit with the period $T=31.7689$, eigenvalues $\lambda_s=-3, \ \lambda_u=1$, nontrivial Floquet exponent $k=-1.7579$, and the eigenvector $v = 0.0006 \hat x +   0.9999 \hat y$. The time series, periodic orbit, and the box $B$ are shown in Fig. \ref{simple_homoclinic}.

\begin{figure*}[!t]
\begin{center}
\includegraphics[width=\textwidth]{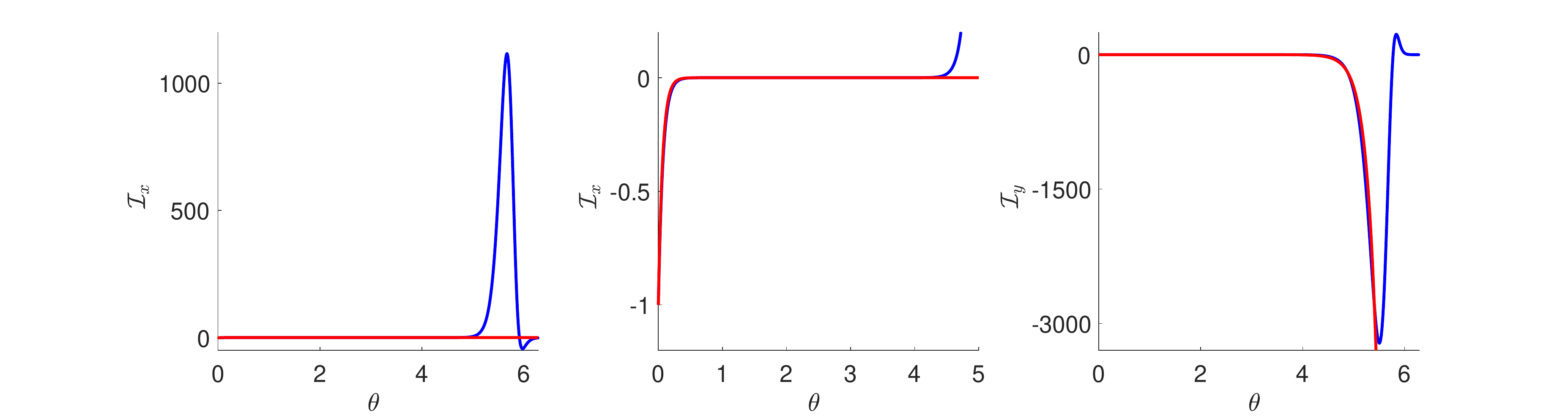}
\end{center}
\caption{IRC for periodic orbit near a homoclinic bifurcation. The left and the right panels show the $x$ and $y$ component of the IRC respectively, the middle panel shows the zoomed in plot of the left panel. The blue line shows the numerically computed IRC, while the red line shows an exponential curve with rate constant given by (\ref{homo_irc}).}
\label{homoclinic_irc}
\end{figure*}

With $\Delta=0.0201$, the trajectory spends 86.5 \% of its period in the box $B$. Fig. \ref{homoclinic_irc} compares the numerically computed IRC with the exponential curve having rate constants from the analytical IRC (\ref{homo_irc}). We see that the numerically computed IRC agrees well with the analytical one in the beginning (inside box $B$), but diverges after. It oscillates quickly back to its initial value at the end of its period, as is expected.

\subsection{Relaxation oscillator}\label{rel_osc}
In a relaxation oscillator, at least one variable evolves at a much faster rate than the other variables. Such oscillators are ubiquitous in conductance-based models of cells, where the gating variables evolve at a much slower rate than the cell membrane potential. A two-dimensional relaxation oscillator can be written as
\begin{eqnarray}
\mu \dot x &=& f(x,y),\qquad 0<\mu\ll1,\\
\dot y &=& g(x,y).
\end{eqnarray}
In the relaxation limit ($\mu \rightarrow 0$), the PRC is given as \cite{izhikevich}
\begin{eqnarray}
\mathcal{Z}(\theta) = -\frac{\omega g_x}{f_x g}\hat x+\frac{\omega }{ g}\hat y.
\end{eqnarray}
Here the functions $g, g_x,$ and $f_x$ are evaluated on the periodic orbit, and thus are functions of $\theta$.
The eigenvector $v$ in the direction of the isochron is given as 
\begin{eqnarray}
v= \frac{-\hat x-\frac{ g_x}{f_x}\hat y}{\sqrt{1+\frac{g_x^2}{f_x^2}}}.
\end{eqnarray}
For computing the adjoint equation for IRC in relaxation limit, we do the following analysis in the spirit of \cite{brow04}. 

Consider an infinitesimal perturbation ${\bf \Delta x}=(\Delta x, \Delta y)$ to the periodic trajectory ${\bf x} \in \gamma(t)$. Then the perturbed trajectory evolves as
\begin{eqnarray}
\mu \dot{\Delta x} &=&f_x \Delta x + f_y \Delta y,\\
\dot {\Delta y} &=&g_x \Delta x +g_y \Delta y.
\end{eqnarray}
This can be written as $A\dot {{\bf \Delta x}}=DF{\bf \Delta x}$, where $A=\begin{bmatrix}\mu & 0\\0 & 1\end{bmatrix}$, and $DF$ is the Jacobian evaluated on the periodic orbit. The isostable shift $\Delta \psi$ by a perturbation $A {{\bf \Delta x}}$  is given by $\Delta \psi = \langle\nabla \psi,A {{\bf \Delta x}}\rangle$, where $\langle\cdot,\cdot\rangle$ is the Euclidean inner product. Its time evolution can be written as
\begin{eqnarray*}
\dot{\Delta \psi}&=&\langle\nabla \dot \psi,A {{\bf \Delta x}}\rangle+\langle\nabla \psi,A\dot {{\bf \Delta x}}\rangle=k\Delta \psi,\\
&\Rightarrow& \langle A^T\nabla \dot \psi,{{\bf \Delta x}}\rangle=\langle kA^T\nabla \psi, {{\bf \Delta x}}\rangle-\langle\nabla \psi,DF{\bf \Delta x}\rangle.
\end{eqnarray*}
This can be written as
\begin{eqnarray}
\mu \mathcal{\dot I}_x&=&(k\mu-f_x) \mathcal{I}_x -g_x \mathcal{I}_y,\\
  \mathcal{\dot I}_y&=&-f_y \mathcal{I}_x +(k-g_y) \mathcal{I}_y,
\end{eqnarray}
where $ \mathcal{I}_x=\partial \psi/ \partial x$, and $ \mathcal{I}_y=\partial \psi/ \partial y$. From the mean of the divergence of the vector field along periodic trajectory, we get the  nontrivial Floquet exponent and multiplier as
\begin{eqnarray}
\lambda&=&\exp{\left(\int_0^T\left(f_x/\mu +g_y\right)dt\right)},\label{floq}\\
k&=&a/\mu +b\label{floq_exponent},
\end{eqnarray}
where $a=\frac{\int_0^Tf_xdt}{T}$, and $b=\frac{\int_0^Tg_ydt}{T}$. We must have $k<0$ for a stable periodic orbit.  This implies that $a< 0$, because otherwise, $k$ would get positive as $\mu \rightarrow 0$. Thus in the relaxation limit, $k \rightarrow -\infty $ and $\lambda \rightarrow 0$, and any perturbation from the periodic orbit gets nullified instantly by the vector field. The adjoint equation for the IRC becomes
\begin{eqnarray}
\mu \mathcal{\dot I}_x&=&(a+\mu b-f_x) \mathcal{I}_x -g_x \mathcal{I}_y,\\
  \mathcal{\dot I}_y&=&-f_y \mathcal{I}_x +(a/\mu+b-g_y) \mathcal{I}_y.
\end{eqnarray}
\begin{eqnarray}
\Rightarrow \mathcal{I}_x&=& \frac{g_x}{a+\mu b-f_x} \mathcal{I}_y+{\cal O}(\mu),\label{irc_relax_x}\\
  \Rightarrow \mu \mathcal{\dot I}_y&=&\left(a+\mu b-\mu g_y-\frac{\mu g_xf_y}{a+\mu b-f_x}\right) \mathcal{I}_y+{\cal O}(\mu^2).\nonumber
\end{eqnarray}
In the relaxation limit ($\mu \rightarrow 0$), we get
\begin{equation}
\left(a-f_x\right)\mathcal{I}_y=0.\label{irc_relax}
\end{equation}
We know from the mean value theorem that there is at least one phase $\theta _i$ where $a= f_x$. Thus the coefficient of $\mathcal{I}_y$ in (\ref{irc_relax}) is nonzero except at $\theta_i$. Thus in order to satisfy the eq. (\ref{irc_relax}), $\mathcal{I}_y$ has to be zero for all $\theta$ except at $\theta_i$ where it can be non-zero. The same can be said about $\mathcal{I}_x$ from eq. (\ref{irc_relax_x}). Thus we can write the the IRC as
\begin{eqnarray}
 \mathcal{I}_{x,y} &=& \left(\underset{i}{\Sigma} \ \mathcal{I}_x(\theta_i)\right) \ \hat x + \left(\underset{i}{\Sigma} \ \mathcal{I}_y(\theta_i)\right) \ \hat y.\label{relaxation_irc}
\end{eqnarray}
It makes sense intuitively that the IRC is zero everywhere except at few points because the periodic orbit is very strongly stable in the relaxation limit (the nontrivial Floquet multiplier is close to zero). Therefore, a perturbation from the periodic orbit gets nullified instantaneously by the stabilizing vector field. This renders the isostable coordinate zero near the periodic orbit, and its gradient zero almost everywhere on the periodic orbit.
\subsubsection{van der Pol oscillator}
\begin{figure*}[!t]
\begin{center}
\includegraphics[width=\textwidth]{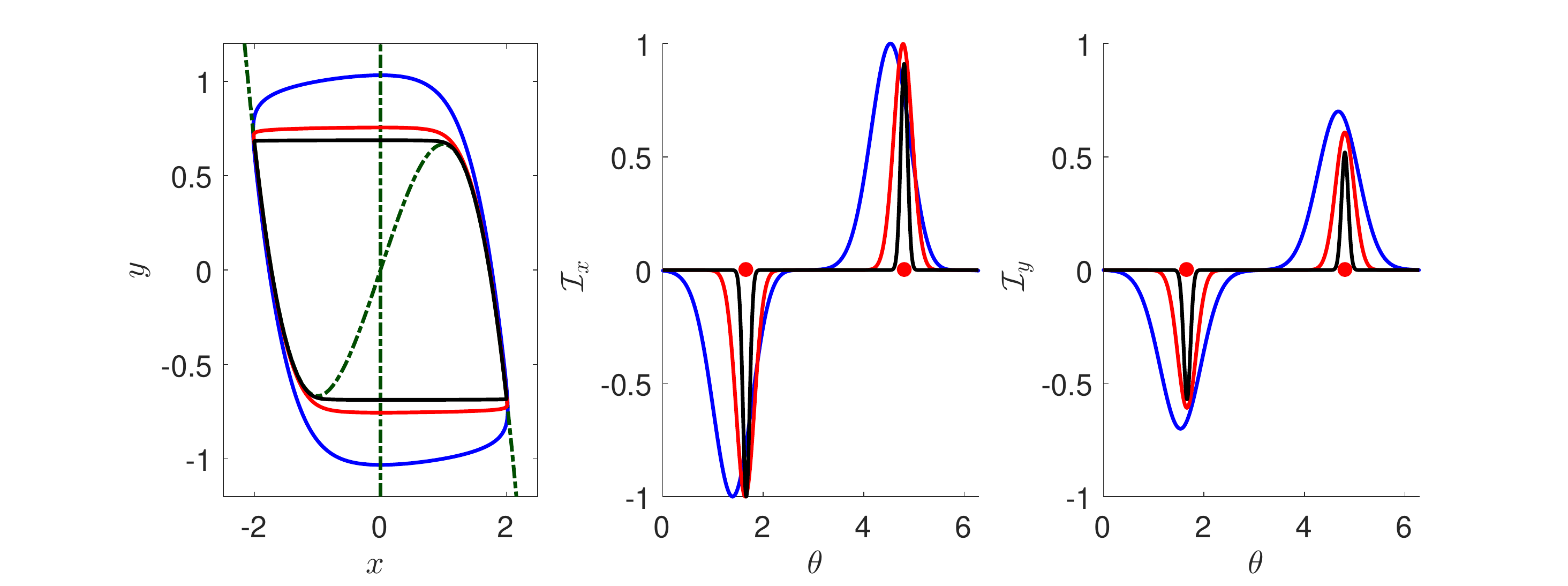}
\end{center}
\caption{van der Pol Oscillator: The left panel plots the periodic orbits and nullclines. The middle (resp., the right) panel plots $\mathcal{I}_x$ (resp., $\mathcal{I}_y$). In all plots, the blue, red and black lines correspond to $\mu=0.1,\ 0.01,$, and $0.001$, respectively. The two red dots in the middle and right panels mark the phases $\theta_1$ and $\theta_2$.}
\label{relaxation}
\end{figure*}

An example of a relaxation oscillator is the van der Pol oscillator \cite{van1920,van_scholar} which can be written as
\begin{eqnarray}
\mu \dot x&=&-y+x-x^3/3,\qquad 0<\mu \ll 1,\\
\dot y &=& x.
\end{eqnarray}
In the relaxation limit ($\mu \rightarrow 0$), we find numerically that $a-f_x$ crosses zero at $\theta_1=1.6567$ and $\theta_2=4.7983$. Thus we expect the IRC to be zero everywhere except these two $\theta_i$ values. We compute periodic orbits and their IRCs for three different values of the parameter $\mu: 0.1, 0.01$, and $0.001$, as shown in Fig. \ref{relaxation}. 
We see from Fig. \ref{relaxation} that as $\mu$ approaches the relaxation limit, the IRC becomes zero everywhere except near the phases $\theta_1$, and $\theta_2$, thus validating our analytical results. Since the IRC is zero everywhere except near 2 points, we do not use the normalization condition of Section \ref{APR}; instead we normalize the IRC by the maximum absolute value of $\{\mathcal{I}_x(\theta_i),\mathcal{I}_y(\theta_i)\}$.

\section{Discussion and Conclusions} \label{section_4}
Standard phase reduction is a crucial tool in the analysis and control of oscillators. It reduces the dimensionality of a system, and can make its control experimentally amenable. However it only allows a small perturbation without the risk of driving the oscillator away from the periodic orbit. This limitation makes it unsuitable for some control purposes, especially when a significant control stimulus is required or when a nontrivial Floquet exponent of the periodic orbit has small magnitude. This necessitates the use of the augmented phase reduction. 

In this article, we have derived expressions for the augmented phase reduction for two distinct systems with a periodic orbit. We find that for a system near homoclinic bifurcation, the IRC is exponential for a large part of its phase. For a relaxation oscillator, the IRC is zero everywhere except at a few points. We simulated dynamic models which are examples of these two systems, and found that their numerically computed IRCs match with their analytical counterparts very closely. From our calculations in \cite{tutorial}, we found that the $\lambda \ - \ \omega$, Hopf, and Bautin normal form systems have sinusoidal PRCs and IRCs. For a model near a SNIPER bifurcation, the PRC never changes sign, while the IRC looks like a skewed sinusoid.

For a strongly stable system, the nontrivial Floquet exponent $k$ goes to $-\infty$. This is the case for relaxation oscillator in the relaxation limit. Thus, any perturbation to the periodic orbit gets nullified instantly. In such a case, it is not necessary to use the augmented phase reduction, instead the standard phase reduction would suffice. On the other hand, for the other five systems, it is better to use the augmented phase reduction over the standard phase reduction, especially when $k$ is a negative number that is small in magnitude.

Table \ref{summary} summarizes the analytical expressions for augmented phase reduction derived in Section \ref{section_3} and in \cite{tutorial}.
      
      \afterpage{%
    \clearpage
    \begin{landscape}
\begin{table}[!t]
\caption{Summary of analytical expressions of augmented phase reduction for six dynamically distinct systems}
\label{summary}
\renewcommand{\arraystretch}{1.5}
 \centering 
\begin{tabular}[c]{|c|c|c|c|}
  \hline
Dynamic Model& PRC &IRC & $k$\\
\hline
$\begin{aligned}[t] \lambda \ - \ \omega \ \  \\ \dot r=G(r),\\ \dot \phi=H(r).\end{aligned}$& $\begin{aligned}[t] \left(-\frac{H'(r_{po})}{G'(r_{po})}\cos\theta-\frac{\sin\theta}{r_{po}}\right)\hat x\\+\left(-\frac{H'(r_{po})}{G'(r_{po})}\sin\theta+\frac{\cos\theta}{r_{po}}\right)\hat y\end{aligned}$ & $\begin{aligned}[t]-\sqrt{1+\frac{{r_{po}}^2H'(r_{po})^2}{G'(r_{po})^2}}\cos\theta~\hat x\\-\sqrt{1+\frac{{r_{po}}^2H'(r_{po})^2}{G'(r_{po})^2}}\sin\theta~\hat y\end{aligned}$ & $G'(r_{po})$\\ 
\hline

$\begin{aligned}[t] \text{Hopf} \ \ \ \ \ \\ \dot r=ar+cr^3,\\ \dot \phi=b+dr^2.\end{aligned}$& $\begin{aligned}[t] \left(\frac{d}{\sqrt{-ac}}\cos\theta+\frac{c}{\sqrt{-ac}}\sin\theta\right)\hat x\\+\left(\frac{d}{\sqrt{-ac}}\sin \theta-\frac{c}{\sqrt{-ac}}\cos\theta\right)\hat y\end{aligned}$ & $\begin{aligned}[t]-\sqrt{1+\frac{d^2}{c^2}}\cos\theta~\hat x\\-\sqrt{1+\frac{d^2}{c^2}}\sin\theta~\hat y\end{aligned}$ & $-2a$\\ 
      
      \hline
      
$\begin{aligned}[t] \text{Bautin} \ \ \ \ \ \ \ \ \\ \dot r=ar+cr^3+fr^5,\\ \dot \phi=b+dr^2+gr^4.\end{aligned}$& $\begin{aligned}[t] \left(-\frac{2dr_{po}+4gr_{po}^3}{a+3cr_{po}^2+5fr_{po}^4}\cos\theta-\frac{\sin\theta}{r_{po}}\right)~\hat x\\+\left(-\frac{2dr_{po}+4gr_{po}^3}{a+3cr_{po}^2+5fr_{po}^4}\sin\theta+\frac{\cos\theta}{r_{po}}\right)~\hat y\end{aligned}$ & $\begin{aligned}[t]-\sqrt{1+r_{po}^2\left(\frac{2dr_{po}+4gr_{po}^3}{a+3cr_{po}^2+5fr_{po}^4}\right)^2}\cos\theta~\hat x\\-\sqrt{1+r_{po}^2\left(\frac{2dr_{po}+4gr_{po}^3}{a+3cr_{po}^2+5fr_{po}^4}\right)^2}\sin\theta~\hat y\end{aligned}$ & $\begin{aligned}[t] a+3cr_{po}^2\\+5fr_{po}^4\end{aligned}$\\
      
      \hline

$\begin{aligned}[t] \text{SNIPER} \ \ \ \\ \dot r = \rho r -r^3,\\ \dot \phi=\eta-\sin\phi.\end{aligned}$& $\begin{aligned}[t] \frac{\cos \theta +\sqrt{\eta^2-1}\sin \theta-1}{\sqrt{\rho}\sqrt{\eta^2-1}}~\hat x\\+\frac{\sin \theta -\sqrt{\eta^2-1}\cos \theta}{\sqrt{\rho}\eta}~\hat y\end{aligned}$ & $\begin{aligned}[t]\frac{\sqrt{\eta^2-1}\sin\theta-(\eta^2-1)\cos\theta}{\sqrt{\eta^2-1}\sin\theta+\cos\theta-\eta^2}~\hat x\\+ \frac{\eta\left(1-\sqrt{\eta^2-1}\sin\theta-\cos\theta\right)}{\sqrt{\eta^2-1}\sin\theta+\cos\theta-\eta^2}~\hat y\end{aligned}$ & $-2\rho$\\ 

 \hline
      
$\begin{aligned}[t] \text{Homoclinic}\\ \dot x = \lambda _u x,\\ \dot y = \lambda_s y.\end{aligned}$& $\begin{aligned}[t]  \mathcal{Z}_{x_0}e^{-\frac{\lambda_u \theta}{\omega}}\hat x+\mathcal{Z}_{y_0}e^{-\frac{\lambda_s \theta}{\omega}}\hat y \end{aligned}$ & $\begin{aligned}[t]\mathcal{I}_{x_0} e^{\frac{\lambda_s \theta}{\omega}}~\hat x + e^{\frac{\lambda_u \theta}{\omega}}~\hat y\end{aligned}$ & $\lambda_s + \lambda_u$\\ 

      \hline
 
$\begin{aligned}[t] \text{Relaxation} \ \ \\ \mu \dot x = f(x,y),\\ \dot y = g(x,y).\end{aligned}$& $\begin{aligned}[c] -\frac{\omega g_x}{f_x g}~\hat x+\frac{\omega }{ g}~\hat y\end{aligned}$ & $\begin{aligned}[c]\left(\underset{i}{\Sigma} \ \mathcal{I}_x(\theta_i)\right) ~\hat x + \left(\underset{i}{\Sigma} \ \mathcal{I}_y(\theta_i)\right)~\hat y\end{aligned}$ & $-\infty$\\ 

      \hline
      \end{tabular}
      \end{table}
\end{landscape}
\clearpage
}

\section*{Acknowledgment}

This work was supported by National Science Foundation Grant No. NSF-1264535/1631170.





\bibliographystyle{siamplain}
\bibliography{ms}

\end{document}